\documentclass[10pt]{amsart}

\usepackage{amsfonts,amsmath,latexsym,amssymb} 
\usepackage[T2A]{fontenc}
\usepackage[utf8]{inputenc}
\usepackage[ukrainian,english]{babel}
\usepackage[dvipsnames]{xcolor}
\usepackage{cite}

\newtheorem{theorem}{Theorem}

\newcommand{\Rd}{{\mathbb R}^d}
\newcommand{\RR}{\mathbb R}

\newcommand{\polarK}{K^\circ}

\newcommand{\KC}{K\cap C}

\allowdisplaybreaks

\begin{document}

\title[\empty]{On some sharp  Landau -- Kolmogorov -- Nagy type inequalities in Sobolev spaces of multivariate functions}

\author{V.~F.~Babenko}
 \address{Department of Mathematical Analysis and Theory of Functions, Oles Honchar Dnipro National University, Dnipro, Ukraine}
\email{babenko.vladislav@gmail.com}

\author{V.~V.~Babenko}
 \address{Department of Mathematics and Computer Science, Drake University, Des Moines, USA}
 \email{vira.babenko@drake.edu}

\author{O.~V.~Kovalenko}
\address{Department of Mathematical Analysis and Theory of Functions, Oles Honchar Dnipro National University, Dnipro, Ukraine}
 \email{olegkovalenko90@gmail.com}

\author{N.~V.~Parfinovych}
\address{Department of Mathematical Analysis and Theory of Functions, Oles Honchar Dnipro National University, Dnipro, Ukraine}
 \email{parfinovich@mmf.dnu.edu.ua}


\begin{abstract}
 For a function $f$ from the Sobolev space $W^{1,p}(C)$ ($C\subset\RR^d$ is an open convex cone), a sharp inequality that estimates $\| f\|_{L_{\infty}}$ via the $L_{p}$-norm of its gradient and a seminorm of the function is obtained.
With the help of this inequality, a sharp inequality is proved, which estimates the ${L_{\infty}}$-norm of the Radon--Nikodym derivative of a charge defined on Lebesgue measurable subsets of 
$C$ via the $L_p$-norm of the gradient of this derivative and a seminorm of the charge. In the case, when $C=\RR_+^m\times \RR^{d-m}$, $0\le m\le d$, we obtain inequalities that estimate the ${L_{\infty}}$-norm of a mixed derivative of a function $f\colon C\to \RR$ using its ${L_{\infty}}$-norm and the $L_p$-norm of the gradient of the function's mixed derivative. 
\end{abstract}
\maketitle

\section{Introduction}
Inequalities for norms of intermediate derivatives of univariate and multivariate functions have an important role in many branches of Analysis and its applications. The main attention is focused on sharp inequalities of this kind, and the results by Landau~\cite{Landau}, Kolmogorov~\cite{Kolmogorov1939} and Nagy~\cite{Nagy} are among the brightest ones in the area. A survey of results and further references can be found in~\cite{BKKP, Babenko22}. 
Inequalities for Radon--Nikodym derivatives of charges defined on Lebesgue measurable subsets of an open cone  $C\subset \RR^d$, as well as for charges defined on measurable subsets of a metric space with a measure were obtained in~\cite{Babenko23a,Babenko23b}.

In this article for a function $f$ from the Sobolev space  $W^{1,p}(C)$ ($C\subset\RR^d$ is an open cone) we obtain a sharp Nagy type inequality that estimates  $\| f\|_{L_{\infty}(C)}$ via the  $L_{p}(C)$-norm of its gradient and some seminorm of the function. Using this inequality, we prove a sharp Landau--Kolmogorov type inequality that estimates the ${L_{\infty}(C)}$-norm of the  Radon--Nikodym derivative via the $L_p(C)$-norm of the gradient of this derivative and the value of some seminorm of the charge. In the case  $C=\RR_+^m\times \RR^{d-m}$, $0\le m\le d$, we obtain inequalities that estimate  $ {L_{\infty}(C)}$-norm of a mixed derivative of a function  $f\colon C\to \RR$ via $ {L_{\infty}(C)}$-norm of the function and the $L_p(C)$-norm of the gradient of the mixed derivative. The obtained inequalities are sharp for $m=0,1$. 

The results of the article together with the known general facts, allow to solve the problem of approximation of an unbounded operator by bounded ones, as well as several related problems (see e.g.~\cite{Babenko23a,Babenko23b}).

\section{Notations and some auxiliary results}

For $x,y\in\Rd$, $d\geq 1$, we denote by $(x,y)$ the dot products of the vectors $x$ and $y$. Let  $K\subset \Rd$ be a bounded open convex symmetric with respect to the origin  $\theta$ set such that $\theta\in{\rm int\,} K$. For $x\in \Rd$ we denote by
$|x|_K$ the norm of the vector $x$ that is generated by the set $K$ i.e., 
$|x|_K=\inf\{\lambda>0\colon x\in \lambda K\}.$ Let also 
$
|z|_{\polarK} = \sup \{(x,z)\colon |x|_K\leq 1\}
$ be the dual norm. Everywhere below $C\subset \RR^d$ is an open convex cone, $\mu$ is the Lebesgue measure in $C$, and for $p\geq 1$, we set  $p' = {p}/(p-1)$.

For a measurable set $Q\subset C$ by $L_p(Q), 1\le p\le \infty,$ we denote the set of measurable functions $f\colon Q\to\RR$ with corresponding norms $\| \cdot\|_{L_p(Q)}$; by  $L_{{\rm loc}}(C)$ the space of locally integrable functions i.e., functions $f\colon C\to\RR$ that are integrable on each compact $Q\subset C$. In the space $L_{\rm loc}(C)$ consider a family of seminorms
\[
\rfloor f\lceil_h
=
\sup\limits_{x\in C}\left|\;\int\limits_{hK\cap C}f(x+u)du\right|, h > 0, \text { and } \rfloor f\lceil =\sup\limits_{h>0}\rfloor f\lceil_h.
\]
By $L_{\rfloor \cdot\lceil_h}(C)$ ($L_{\rfloor \cdot\lceil}(C)$) we denote the space of functions  $f\in L_{\rm loc}(C)$ such that $\rfloor f\lceil_h<\infty$ (resp. $\rfloor f\lceil<\infty$). It is clear that  $L_1(C)\subset L_{\rfloor \cdot\lceil}(C)$.

For a locally integrable function $f$ by  $\nabla f$ we denote its gradient, where the derivatives are understood in the distributional sense. Let  $Q\subset \Rd$ be an open set. By $W^{1,p}(Q), 1\le p\le\infty,$  we denote the Sobolev space of functions $f\colon Q\to \RR$ such that all their first order derivatives belong to $L_p(Q)$.
For  $f\in W^{1,p}(Q)$ one has $|\nabla f|_{\polarK} \in L_p(Q)$.

For  $h>0$ we consider the function $g_h\colon (0,h)\to\RR$,
\[
g_h(u)=\frac 1{d\cdot \mu(K\cap C)}\left(\frac 1{u^{d-1}}-\frac u{h^d}\right).
\]
Using~\cite[Lemma~3]{Babenko22},  and transforming the obtained integral we obtain
$$\|g_h(|\cdot|_K)\|_{L_{p'}(hK\cap C)} = (d\cdot\mu(K\cap C))^{1/{p'}}\left(\int_0^ht^{d-1}g_h^{p'}(t)dt\right)^{1/{p'}}=
 \frac{A h^{1-\frac dp}}{\mu^{-\frac 1{p}}(\KC)},$$ 
where (below $B(\cdot,\cdot)$ denotes the Euler $B$-function)
\begin{equation}\label{constant}
A = A(d,p) = d^{-1}  B^{\frac 1{p'}}\left(1-\frac{(d-1)p'}{d}, p'+1\right).
\end{equation} 
Define an operator $S_h\colon L_{\rfloor \cdot\lceil_h}(C)\to L_\infty(C)$,
\[
S_hf(x)=\frac 1{h^d\mu(\KC)}\int_{{h\KC}}f(x+y)dy.
\]
Now it becomes obvious that the following theorem is a partial case of~\cite[Theorem~2]{Babenko22}.
\begin{theorem}\label{th::knownOstrowski}
Let $p\in (d,\infty]$, $h>0$ and  $f\in W^{1,p}(h \KC)$. Then
\begin{equation*}
 \left| f(\theta)-S_hf(\theta)
\right|\leq 
 \left\|g_{h}\left( |\cdot|_K\right) \right\|_{L_{p'}(h\KC)}
\||\nabla f|_{\polarK}\|_{L_p(h\KC)}.
\end{equation*}
The inequality is sharp. It becomes equality for the function $\alpha\cdot f +\beta$, where $\alpha,\beta\in\RR$ and
$$
f(y)= \int_0^{|y|_{K}}g_{h}^{p'-1}(u)du, \, y\in h\KC.
$$
\end{theorem}

\section{Nagy type inequality}
Using Theorem~\ref{th::knownOstrowski}, for all  $x\in C$ one obtains
\begin{multline}\label{knownEstimate}
   \left|f(x) - S_hf(x)\right|
   \leq 
   \left\|g_{h}(|\cdot|_K) \right\|_{L_{p'}(hK\cap C)} \|| \nabla f(x+\cdot)|_{\polarK}\|_{L_p( h\KC)}
\\\le 
\left\|g_{h}(|\cdot|_K) \right\|_{L_{p'}(hK\cap C)}  
\||\nabla f|_{\polarK}\|_{L_p(C)}.
\end{multline}

\begin{theorem}\label{th::NagyInequality}
If $h> 0$, $p\in (d,\infty]$ and $f\in W^{1,p}(C)\cap  L_{\rfloor \cdot\lceil_h}(C)$, then $f\in L_\infty(C)$ and the following inequality holds:
\begin{equation}\label{mainInequality}
\|f\|_{L_\infty(C)}
 \le\left\|g_{h}(|\cdot|_K) \right\|_{L_{p'}(hK\cap C)} \||\nabla  f|_{\polarK}\|_{L_p(C)}+
\mu^{-1}(K\cap C)h^{-d}\rfloor f\lceil_h.
\end{equation}
Inequality~\eqref{mainInequality} is sharp. It becomes equality on the function
\begin{equation}\label{extremalFunction}
	f_{e,h}(y)= \begin{cases} \int_{|y|_{K}}^{h}g_{h}^{p'-1}(u)du , & y\in h\KC,\\
    0, & y\in C\setminus hK.
 \end{cases}
\end{equation}
For each $f\in W^{1,p}(C)\cap  L_{\rfloor \cdot\lceil}(C)$  the following multiplicative inequality holds:
\begin{equation}\label{multiplicativeNagy}
\|f\|_{L_\infty(C)}
\leq 
a(d,p) \mu^{-\frac {\alpha}{d}}(\KC)\rfloor f\lceil ^{1-\alpha} \||\nabla f|_{\polarK}\|_{L_p(C)} ^{\alpha}, 
\end{equation}
where
\begin{equation}\label{constantsDefinition}
    \alpha = \frac{pd}{p+(p-1)d}, a(d,p) = \left(\frac{(p-d)A(d,p)}{pd}\right)^{\alpha}\left(\frac{pd}{p-d}+1\right), 
\end{equation}
and $A(d,p)$ is defined in~\eqref{constant}.
Inequality~\eqref{multiplicativeNagy} is sharp. It becomes equality for each function $f_{e,h}$, $h>0$.

If $f\in W^{1,p}(C)\cap  L_{1}(C)$, then in inequalities~\eqref{mainInequality} and~\eqref{multiplicativeNagy}, $\rfloor f\lceil_h$ and $\rfloor f\lceil$ respectively can be substituted by $\| f\|_{L_1(C)}$, and the obtained inequalities remain sharp.
\end{theorem}

\begin{proof}
Using~\eqref{knownEstimate}, for each $x\in C$ we have 
\begin{multline}\label{trngInequality}
   |f(x)| \leq \left|f(x) - S_hf(x)\right| + \mu^{-1}(K\cap C)h^{-d}\left|\int_{h\KC} f(x+y)dy\right|
    \\\leq 
   \left\|g_{h}(|\cdot|_K) \right\|_{L_{p'}(hK\cap C)} \||\nabla  f|_{\polarK}\|_{L_p(C)}+
\mu^{-1}(K\cap C)h^{-d}\rfloor f\lceil_h,
\end{multline}
which implies inequality~\eqref{mainInequality}. Next we prove its sharpness. Since the function $g_{h}$ is non-negative, one has 
$
   |x|_K\leq |y|_K\implies f_{e,h}(x)\geq f_{e,h}(y).
$
The function $f_{e,h}$ is continuous. Hence $\|f_{e,h}\|_{L_\infty(C)} = f(\theta)$ and, moreover, the first inequality in~\eqref{trngInequality} becomes equality for $x =\theta$ and the function $f_{e,h}$. According to Theorem~\ref{th::knownOstrowski}, the first inequality in~\eqref{knownEstimate} for $x = \theta$ becomes equality for the function $f_{e,h}$. The function $f_{e,h}$ vanishes outside the set $h\KC$, and hence
$$
\||\nabla f_{e,h}|_{\polarK}\|_{L_p(h\KC)} = \||\nabla f_{e,h}|_{\polarK}\|_{L_p(C)} \text{ and } \rfloor f_{e,h}\lceil_h = \int_{h\KC} f_{e,h}(y)dy.
$$
Thus the second inequalities in~\eqref{knownEstimate} and~\eqref{trngInequality} also become equalities, and hence inequality~\eqref{mainInequality} becomes equality on the function  $f_{e,h}$.

From~\eqref{constant} and~\eqref{mainInequality} it follows that for all  $f\in W^{1,p}(C)\cap  L_{\rfloor \cdot\lceil}(C)$ and $h>0$
\begin{equation}\label{additiveNagy}
\|f\|_{L_\infty(C)}
\le A(d,p) \mu^{-\frac 1{p}}(\KC) h^{1 - \frac{d}{p}}
\||\nabla f|_{\polarK}\|_{L_p(C)} 
+
\mu^{-1}(K\cap C)h^{-d}\rfloor f\lceil.
\end{equation}
 Moreover, for all $h>0$, $\rfloor f_{e,h}\lceil_h = \rfloor f_{e,h}\lceil$, hence inequality~\eqref{additiveNagy} is sharp and becomes equality for the function $f_{e,h}$. 

In order to prove~\eqref{multiplicativeNagy} it is sufficient to plug 
$$
h 
=
\left(\frac{pd\mu^{-\frac 1{p'}}(\KC) \rfloor f\lceil}{(p-d) A(d,p)\||\nabla f|_{\polarK}\|_{L_p(C)}}\right)^{\frac{p'}{p'+d}}
$$
into the right-hand side of~\eqref{additiveNagy}. 

Next we prove sharpness of inequality~\eqref{multiplicativeNagy}. Indeed, from this inequality, using Young's inequality it is not hard to obtain that for each $h>0$ inequality~\eqref{additiveNagy} holds. It turns into equality for the function $f_{e,h}$. Hence inequality~\eqref{multiplicativeNagy} also becomes equality on the function $f_{e,h}$. 

The last statement of the theorem is obvious.
\end{proof}

\section{Landau--Kolmogorov type inequalities for charges}
By $\mathfrak{N}(C)$ we denote the family of charges $\nu$ defined on Lebesgue measurable subsets of the set $C$ that are absolutely continuous with respect to the Lebesgue measure  $\mu$, see e.g.~\cite[Chapter~5]{Berezanski}. 
The Radon--Nikodym derivative of a charge $\nu$ with respect to the measure $\mu$ will be denoted by $D_\mu\nu$.
The family $\mathfrak{N}(C)$ is a linear space with respect to the standard addition and multiplication by a real number. For $h>0$ and $\nu\in \mathfrak{N}(C)$ we consider seminorms
$\rceil \nu\lfloor_{h}=\sup_{x\in C}|\nu(x +hK)|$ and $\rceil \nu\lfloor=\sup_{h>0}\rceil \nu\lfloor_{h}$.
It is clear that
$\rceil \nu\lfloor_h=\rfloor D_\mu\nu\lceil_h$ i $\rceil \nu\lfloor=\rfloor D_\mu\nu\lceil.$
By $\mathfrak{N}_{\rceil \cdot\lfloor_h}(C)$ ($\mathfrak{N}_{\rceil \cdot\lfloor}(C)$) we denote the set of charges $\nu\in\mathfrak{N}(C)$ that have finite seminorm $\rceil \cdot\lfloor_h$ (resp. $\rceil \cdot\lfloor$). 

Applying Theorem~\ref{th::NagyInequality} to the function $f=D_\mu\nu$, we obtain the following result.
\begin{theorem}\label{th::LKForCharges}
If $h>0$ and the charge $\nu\in \mathfrak{N}_{\rceil \cdot\lfloor_h}(C)$ is such that $D_\mu\nu\in W^{1,p}(C)$, then
\begin{equation}\label{mainInequalityCharges}
\|D_\mu\nu\|_{L_\infty(C)}
\leq   
A \mu^{-\frac 1{p}}(\KC) h^{1 - \frac{d}{p}}
\||\nabla D_\mu\nu|_{\polarK}\|_{L_p(C)} +
\mu^{-1}(K\cap C)h^{-d}\rceil \nu\lfloor_h,
\end{equation}
where $A(d,p)$ is defined in~\eqref{constant}.
Inequality~\eqref{mainInequalityCharges} is sharp. It becomes equality for the charge $\nu_{e,h}$ such that $D_\mu\nu_{e,h}=f_{e,h}$, where the function $f_{e,h}$ is defined in~\eqref{extremalFunction}.
 If the charge $\nu\in \mathfrak{N}_{\rceil \cdot\lfloor}(C)$ is such that $D_\mu\nu\in W^{1,p}(C)$, then the following multiplicative inequality holds.
\begin{equation}\label{multiplicativeNagyCharges}
\| D_\mu\nu\|_{L_\infty(C)}
\leq 
a(d,p) \mu^{-\frac {\alpha}{d}}(\KC)\rfloor \nu\lceil ^{1-\alpha} \||\nabla D_\mu\nu|_{\polarK}\|_{L_p(C)} ^{\alpha}, 
\end{equation}
where $\alpha$ and $a(d,p)$ are defined in~\eqref{constantsDefinition}. Inequality~\eqref{multiplicativeNagyCharges} is sharp. It becomes equality for each charge $\nu$ such that $D_\mu\nu=f_{e,h}$, $h>0$.
\end{theorem}

\section{Inequalities that contain the $L_p$-norm of the gradient of the mixed derivative}
Let $1\le m\le d$, $C=\RR^d_{m,+}=\RR^m_+\times \RR^{d-m}$, and $K=(-1,1)^d$. Then $hK\cap C=(0,h)^m \times (-h,h)^{d-m}$, $h>0$, and $|x|_K=\max\{|x_1|,\ldots,|x_d|\}$. 
For ${\bf I}=(1,\ldots,1)\in \RR^d$, a locally integrable function $f\colon C\to\RR$ and the standard basis  $\{e_i\}_{i=1}^d$ in $\RR^d$ set
$\partial_{\bf I}f=\frac{\partial ^d f}{\partial x_1\ldots\partial x_d}$
(the derivatives are understood in the distributional sense),
\[
\Delta^+_{i,h}f(x)=f(x+he_i)-f(x)
\text{ and }
\Delta_{i,h}f(x)=f(x+he_i)-f(x-he_i).
\]
According to the Fubini theorem, for almost all $x\in \RR^d_{m,+}$ we have
\begin{equation}\label{mixdiff}
\int_{x+hK\cap C} \partial_{\bf I}f(u)du=(\Delta^+_{1,h}\circ\ldots\circ\Delta^+_{m,h}\circ\Delta_{m+1,h}\circ\ldots\circ\Delta_{d,h})f(x). 
\end{equation}
It is easy to see that for the operator  $\mathfrak{S}_{h,m}\colon L_\infty(C)\to L_\infty(C)$
\[
\mathfrak{S}_{h,m}f(x)=\frac 1{2^{d-m}h^d}\left(\Delta^+_{1,h}\circ\ldots\circ\Delta^+_{m, h}\circ\Delta_{m+1,h}\circ\ldots\circ\Delta_{d,h}\right)f(x),
\]
\begin{equation}\label{NormOfS}
\mathfrak{S}_{h,m}f(x)=S_h\partial_{\bf I}f(x)\,\forall f\in L_\infty(C),\text{ and 
 } \|\mathfrak{S}_{h,m}\|_{L_\infty(C)\to L_\infty(C)} = 2^mh^{-d}.
\end{equation}
\begin{theorem}\label{th::mixedDerivatives}
For $h>0$, $K=(-1,1)^d$, $C = \RR^d_{m,+}$ and the function $f\in L_\infty(C)$  such that $\partial_{\bf{I}}f\in W^{1,p}(C)$, the following inequality holds:
\begin{equation}\label{partialLp}
\|\partial_
{\bf{I}}f\|_{L_\infty(C)}
\leq 
A(d,p){h^{1 - \frac{d}{p}}}{2^\frac{m-d}p}
\||\nabla \partial_{\bf I}f|_{\polarK}\|_{L_p(C)} 
+
{2^m}{h^{-d}}\| f\|_{L_\infty(C)},
\end{equation}
where $A(d,p)$ is defined in~\eqref{constant}. It can be rewritten in the following multiplicative form:
\begin{equation}\label{partialLpMultiplicative}
\|\partial_
{\bf{I}}f\|_{L_\infty(C)}
\leq 
a(d,p) 2^{\alpha\left(\frac md - \frac dp\right)}\| f\|_{L_\infty(C)} ^{1-\alpha} \||\nabla  \partial_{\bf I}f |_{\polarK}\|_{L_p(C)} ^{\alpha}, 
\end{equation}
where $\alpha$ and $a(d,p)$ are defined in~\eqref{constantsDefinition}.
For $m=0$ and $m = 1$ inequalities~\eqref{partialLp} and~\eqref{partialLpMultiplicative} are sharp.
\end{theorem}

\begin{proof}
Taking into account that $\mu(\KC)=2^{d-m}$, and formulae~\eqref{knownEstimate} and~\eqref{NormOfS}, one has
\begin{multline*}
\|\partial_{\bf{I}}f\|_{L_\infty(C)}
\leq 
\|\partial_{\bf{I}}f-\mathfrak{S}_{h,m}f\|_{L_\infty(C)} 
+ \|\mathfrak{S}_{h,m}\|_{L_\infty(C)\to L_\infty(C)}\|f\|_{L_\infty(C)}
\\=
\|\partial_{\bf{I}}f-S_h \partial_{\bf{I}}f\|_{L_\infty(C)} 
+ 2^{m}h^{-d}\|f\|_{L_\infty(C)}
\\\leq A(d,p) \cdot 2^{\frac{m-d}p} h^{1-\frac dp}  
\||\nabla \partial_{\bf{I}}f|_{\polarK}\|_{L_p(C)} +  2^{m}h^{-d}\|f\|_{L_\infty(C)},
\end{multline*}
and inequality~\eqref{partialLp} is proved. From inequality~\eqref{multiplicativeNagy} for the function $\partial_{\bf{I}}f$ and equality~\eqref{mixdiff}, we obtain 
\begin{multline*}
\|\partial_{\bf{I}}f\|_{L_\infty(C)}
\leq 
a(d,p) \mu^{-\frac {\alpha}{d}}(\KC)\rfloor \partial_{\bf{I}}f\lceil ^{1-\alpha} \||\nabla \partial_{\bf{I}}f|_{\polarK}\|_{L_p(C)} ^{\alpha}
\\ \leq
a(d,p) 2^{\frac{(m-d)\alpha}{d}}\left(2^d\| f\|_{L_\infty(C)}\right)^{1-\alpha} \||\nabla \partial_{\bf{I}}f|_{\polarK}\|_{L_p(C)} ^{\alpha}
\\=
a(d,p) 2^{\alpha\left(\frac md - \frac dp\right)}\| f\|_{L_\infty(C)} ^{1-\alpha} \||\nabla  \partial_{\bf I}f |_{\polarK}\|_{L_p(C)} ^{\alpha},
\end{multline*}
and inequality~\eqref{partialLpMultiplicative} is proved.
We prove sharpness of inequalities~\eqref{partialLp} and~\eqref{partialLpMultiplicative} for $m = 0$. For the function $f_{e,h}$ defined in~\eqref{extremalFunction}, and the function
$$
F_{e,h}(x) = \int_0^{x_1}\ldots\int_0^{x_d}f_{e,h}(u)du, 
$$
we have $\partial_{\bf I}F_{e,h} = f_{e,h}$, $|\nabla\partial_{\bf I} F_{e,h}(\cdot)|_{\polarK}=|\nabla f_{e,h}(\cdot)|_{\polarK}$, and due to the symmetry considerations,
\begin{equation}\label{FehNorm}
   2^d \|F_{e,h}\|_{L_\infty(C)} = 2^d\int_{(0,h)^d}f_{e,h}(u)du = \int_{(-h,h)^d}f_{e,h}(u)du = 
\rfloor f_{e,h}\lceil_h = 
\rfloor f_{e,h}\lceil.
\end{equation}
Since for the function  $f_{e,h}$ inequalities~\eqref{mainInequality} and~\eqref{multiplicativeNagy} become equalities, taking into account~\eqref{FehNorm}, we obtain sharpness of inequalities~\eqref{partialLp} and~\eqref{partialLpMultiplicative} for $m = 0$.

Next we prove sharpness of inequalities~\eqref{partialLp} and~\eqref{partialLpMultiplicative} in the case $m = 1$. In this case $h\KC = (0,h)\times (-h,h)^{d-1}$ and $\mu (\KC) = 2^{d-1}$. There exists a number $a\in (0,h)$ such that 
\begin{equation}\label{aChoice}
\int_{(0,a)\times (-h,h)^{d-1}} f_{e,h}(u)du = \int_{(a,h)\times (-h,h)^{d-1}} f_{e,h}(u)du 
.   
\end{equation}
Consider the function 
$
G_{e,h}(x) = \int_a^{x_1}\int_0^{x_2}\ldots\int_0^{x_d}f_{e,h}(u)du.
$
We have $\partial_{\bf I}G_{e,h} = f_{e,h}$ and hence $|\nabla\partial_{\bf I} G_{e,h}(\cdot)|_{\polarK}=|\nabla f_{e,h}(\cdot)|_{\polarK}$.

The hyperplanes $x_1 = a$ and $x_j=0$, $j=2,\ldots, d$, split the set $(0,h)\times (-h,h)^{d-1}$ into $2^{d}$ parallelepipeds  $\Pi_1,\ldots, \Pi_{2^d}$; due to the symmetricity with respect to the coordinate hyperplanes of the plot of the function  $f_{e,h}$  and equality~\eqref{aChoice}, we have
$$
\int_{\Pi_i}f_{e,h}(u)du = \frac{1}{2^d}\int_{(0,h)\times (-h,h)^{d-1}} f_{e,h}(u)du, i =1,\ldots, 2^d,
$$
and hence 
$
\|G_{e,h}\|_{L_\infty(C)} =  \frac{1}{2^d}\int_{(0,h)\times (-h,h)^{d-1}} f_{e,h}(u)du =
\frac{1}{2^d}
\rfloor f_{e,h}\lceil_h = 
\frac{1}{2^d}\rfloor f_{e,h}\lceil.
$
From these equalities and the fact that inequality~\eqref{mainInequality} and~\eqref{multiplicativeNagy} become equalities for the function $f_{e,h}$, we obtain sharpness of inequalities~\eqref{partialLp} and~\eqref{partialLpMultiplicative} for $m = 1$.
\end{proof}

\section{Some applications}

Let $X$, $Y$ and $Z$ be linear spaces with seminorms $\|\cdot\|_X$, $\|\cdot\|_Y$ and $\|\cdot\|_Z$ respectively. A linear operator $S\colon X\to Y$ is called bounded (and we write $S\in \mathcal{L}(X,Y)$), if
\[
\| S\|_{X\to Y}=\sup\{\| Sx\|_Y\colon \| x\|_X\le 1\}<\infty.
\]

Let $A\colon X \to Y, B\colon X\to Z$ be two homogeneous operators with domains of definition $D_A, D_B\subset X$, $D_B\subset D_A$. Set $\mathfrak{M}=\{ x\in D_B \colon \| Bx\|_Z\le1\}$.
For the operator $A$ and an operator $S\in \mathcal{L}(X,Y)$  set
$$
    U(A,S; \mathfrak{M}):=\sup\{\|Ax-Sx\|_Y\colon x\in \mathfrak{M}\}.
$$
For arbitrary  $S\in \mathcal{L}(X,Y)$ and each  $x\in D_B$ the following Landau--Kolmogorov--Nagy type inequality holds:
\begin{equation}\label{abstraddineq} 
 \| Ax\|_Y\le  \| Ax-Sx\|_Y+\| S\|_{X\to Y}\|x\|_X
\le U(A,S;\mathfrak{M})\| Bx\|_Z+\|S\|_{X\to Y}\cdot \|x\|_X.
\end{equation}

From Theorem~\ref{th::NagyInequality} it follows that for the operators $A\colon f\mapsto f, B\colon f\mapsto |\nabla f|_{K^\circ}$ and $S\colon f\mapsto S_hf$ inequality~\eqref{abstraddineq}
 becomes equality on the functions $f=f_{e,h}$.
From Theorem~\ref{th::LKForCharges} it follows that for $A\colon \nu\mapsto D_\mu\nu, B\colon \nu\mapsto |\nabla D_\mu\nu|_{K^\circ}$ and $S\colon \nu\mapsto \frac{\nu(x+hK\cap C)}{\mu(hK\cap C}$ inequality~\eqref{abstraddineq}
becomes equality for the charge $\nu=\nu_{e,h}$.
Finally, from Theorem~\ref{th::mixedDerivatives} it follows that for $A\colon f\mapsto \partial_{\bf{I}}f, B\colon f\mapsto |\nabla \partial_{\bf{I}}f|_{K^\circ}$ and $S\colon f\mapsto \mathfrak{S}_{h,m}f$  inequality~\eqref{abstraddineq}
becomes equality for $f=F_{e,h}$ in the case $m=0$, and for $f=G_{e,h}$ in the cases $m=1$.

These observations together with known general facts (see e.g.~\cite[Theorems 1,2]{Babenko23a} and~\cite[Theorem~1]{Babenko23b}),
allow to solve the problems of approximation of the corresponding unbounded operators by bounded ones and related problems.

\bibliographystyle{unsrt}
\bibliography{mybibfile}

\end{document}